Certain Properties of Pythagorean Triangles
involving the interior diameter $2\rho$,
and the exterior diameters $2\rho_\alpha$, $2\rho_\beta$, $2\rho_\gamma$

Part II: The legs case

Konstantine 'Hermes' Zelator
Department of Mathematics
College of Arts and Sciences
Mail Stop 924
The University of Toledo
Toledo, OH 43606-3390
USA

## *1. Introduction*

Given three positive integers's $\alpha, \beta, \gamma$, congruent triangles will exist on the two-dimensional Euclidean plane, with side lengths $\alpha, \beta, \gamma$; as long as the three integers satisfy the triangle inequalities $\alpha + \beta > \gamma, \alpha + \gamma > \beta, \gamma + \beta > \alpha$. And depending on the type of integers $\alpha, \beta, \gamma$ involved, one may be able to construct such a representative triangle (of the said family of congruent triangles) by using ruler and compass; or some other means. Even if such a triangle is not constructible (by using traditional instruments such as a ruler and a compass), it will certainly exist in a mathematical sense. It is also easy to see that there exist infinitely many integer-sided triangles with one side length being a perfect or integer square; and one of the four positive integers $\alpha + \beta + \gamma$ (perimeter), $-\alpha + \beta + \gamma$, $\alpha - \beta + \gamma$, $\alpha + \beta - \gamma$ being also and integer square. Indeed, if we take $\alpha = k^2$ for some positive integer k, and say, $\alpha + \beta + \gamma = \ell^2$, where $\ell$ is also a positive integer. Then, the numbers $\beta$ and $\gamma$ can be chosen so that, $\beta + \gamma = \ell^2 - k^2$ (with $\ell > k) > k^2$, which in turn requires that the positive integers $\ell$ and k satisfy $\ell > k\sqrt{2} \Leftrightarrow \ell^2 > 2k^2$ ), and $-\kappa^2 < \beta - \gamma < \kappa^2$ or equivalently, $|\beta - \gamma| < \kappa^2$. For each value of $t = -\kappa^2 + 1, -k^2, ..., 0, 1, 2, ..., \kappa^2 - 2, \kappa^2 - 1$, a triangle with integer sidelengths

$\alpha = k^2, \beta = \dfrac{\ell^2 - k^2 + t}{2}, \gamma = \dfrac{\ell^2 - k^2 - t}{2}$, is generated; since, as it is easily seen, the three triangle inequalities are satisfied.

The four real numbers $\alpha + \beta + \gamma, -\alpha + \beta + \gamma, \alpha - \beta + \gamma, \alpha + \beta - \gamma$, are also recognizable as being factors in the under the radical quantity in Heron's formula for the area of a triangle.

As we will see in Section 2 of this paper, if we divide twice the area of a triangle by each of the numbers $\alpha + \beta + \gamma, -\alpha + \beta + \gamma, \alpha - \beta + \gamma, \alpha + \beta - \gamma$, we obtain respectively $2\rho, 2\rho_\alpha, 2\rho_\beta, 2\rho_\gamma$. These are the diameters of four important circles: the triangle's inner or inscribed circle and the three exterior circles, each of which is tangential to one side of the triangle and tangential to the two straight lines containing the other two sides (but not tangential to the other two sides themselves). Each of the three centers (of these three circles) is the intersection of one internal bisector (i.e. bisector of one of the triangles internal angle) and two external bisectors (i.e. bisectors of two external triangle angles).

In section 3, we easily show that if a triangle $AB\Gamma$ is a right one with the $90^\circ$ degree angle at the vertex A; then, in fact,

$$2\rho = \beta + \gamma - \alpha, 2\rho_\alpha = \alpha + \beta + \gamma, 2\rho_\beta = \alpha + \beta - \gamma, 2\rho_\gamma = \alpha + \gamma - \beta$$

(so in this case, $\alpha$ is the length of the hypotenuse).

Which brings us to Section 4, which is the focus of this paper, namely, Pythagorean triangles. Therein we obtain simple formulas for $2\rho, 2\rho_\alpha, 2\rho_\beta, 2\rho_\gamma$, in terms of the two integer parameters m and n that generate the family of primitive Pythagorean Triangles.

In Section 5, we state seven well-known results from number theory. Results 3-7 can be found in standard texts and books of number theory we (we offer references [2] and [3] for this); the relevant material covered in the first half of a first course in elementary number theory. In Result 1, we state the parametric formulas which describe the entire family of solutions of the diophantine equation

$x^2+2y^2=z^2$, with (x,y)=1. This is well-known and can be found in E.L. Dickson's landmark book (see [1]). It can also be found in [4].

In Result 2, we state the parametric formulas that describe the entire family of solutions to the diophantine equation $x^2+y^2=2z^2$, with (x,y)=1. More on this in the last section of this paper, Section 8.

In Section 6, we establish Proposition 1: there exist infinitely many primitive Pythagorean triangles with one leg being an integer square. This is a well-known result, but we need to have it at our disposal for what it follows later. This then brings to the central question of this paper. Which is, what can be said about Pythagorean triangles which have one sidelength being an integer square, and one of the four integers $2\rho, 2\rho_\alpha, 2\rho_\beta, 2\rho_\gamma$, also a square? In Section 7, we prove Theorem 1: that in a primitive Pythagorean triangle, with a hypotenuse length $\alpha$; $\beta$ the even leg length, and $\gamma$ the odd leg length; each of the following two combinations is impossible:

1) $\gamma$ a square and $2\rho_\alpha$ also a square.
2) $\gamma$ a square and $2\rho_\gamma$ also a square.

Note that if in a primitive Pythagorean triangle we require that one leg has length an integer square; and that one of the four diameters $2\rho, 2\rho_\alpha, 2\rho_\beta, 2\rho_\gamma$ is also an integer square. Then there are exactly eight combinations. According to the above, two of these combinations are not possible. The remaining combinations are:

3) $\beta$ a square and $2\rho$ also a square.
4) $\beta$ a square and $2\rho_\alpha$ also a square.
5) $\gamma$ a square and $2\rho$ also a square.
6) $\beta$ a square and $2\rho_\beta$ also a square.
7) $\beta$ a square and $2\rho_\gamma$ also a square.
8) $\gamma$ a square and $2\rho_\beta$ also a square.

In Theorem 2, also in Section 7, we show that all of the six combinations are possible, and we parametrically describe each of the six families of such primitive Pythagorean triangles. We also present 12 numerical examples.

We close this paper with Section 8, in which we include some historical comments about the diophanline equation $x^2+y^2=2z^2$, $(x,y)=1$; as well as a brief sketch of a derivation of the general solution of this equation.

*Notation*

1. We will use the standard **notation (a,b)** to denote the greatest common divisor of two integers *a* and *b*.

2. To denote a particular solution to a diophantine equation we will use the bracket notation. For example, $\{x_0, y_0, z_0\}$, can stand for a solution to a three-variable diophantine equation in *x*, *y*, and *z*.

3. The symbol $\mathbb{Z}$ will denote the set of integers; $\mathbb{Z}^+$ the set of positive integers.

4. The notation a|b, will mean that the integer *a* is a divisor of the integer *b*. The same meaning can be conveyed by using the phrases "*a* divides *b*"; or "*b* is divisible by *a*"; or "*a* exactly divides *b*"; or the language of congruences, $b \equiv 0 \pmod{a}$.

5. The notation $\overline{AB}$ will mean "line segment AB", and $|AB|$ will denote the length of that line segment.

## 2. Formulas for the diameters $2\rho, 2\rho_\alpha, 2\rho_\beta, 2\rho_\gamma$ in terms of the sidelengths $\alpha, \beta, \gamma$.

From Figure 1, we easily see that

$$\text{area (of triangle)} \ \overset{\Delta}{ABI} = \frac{1}{2}\rho \bullet \gamma$$

$$\text{area} \ \overset{\Delta}{BI\Gamma} = \frac{1}{2}\rho \bullet \alpha$$

$$\text{area} \ \overset{\Delta}{BI\Gamma} = \frac{1}{2}\rho \bullet \beta$$

$$\Rightarrow \quad \text{area (of triangle)} \ \overset{\Delta}{ABI} = \frac{1}{2}\rho \bullet (\alpha + \beta + \gamma).$$

Also, we have
$$\begin{cases} x + y = \gamma \\ y + z = \alpha \\ z + x = \beta \end{cases} \Longleftrightarrow \text{(algebraically solving)} \quad \begin{cases} x = \frac{\beta + \gamma - \alpha}{2} \\ y = \frac{\alpha + \gamma - \beta}{2} \\ z = \frac{\alpha + \beta - \gamma}{2} \end{cases}$$

Next, $\overline{|AM|} = \overline{|AK|}$ (since $M$ and $K$ are points of tangency of straightlines $AM$ and $AK$, with the circle of center $O_\alpha$); which gives $\gamma + \overline{|BM|} = \beta + \overline{|\Gamma K|}$. But also, $\overline{|BG|} = \overline{|BM|}$ and $\overline{|\Gamma G|} = \overline{|\Gamma K|}$; thus since $\overline{|BG|} + \overline{|\Gamma G|} = \alpha$; it follows that $\overline{|BM|} + \overline{|\Gamma K|} = \alpha$. Combining this with the above $\overline{|BM|} - \overline{|\Gamma K|} = \beta - \gamma$; n easy solving produces $\overline{|BM|} = \frac{\beta + \alpha - \gamma}{2}$ and $\overline{|\Gamma K|} = \frac{\alpha + \gamma - \beta}{2}$

Moreover from right triangle $O_\alpha \overset{\Delta}{MA}$ we see that, $\tan\theta = \frac{\rho_\alpha}{\overline{|AM|}} = \frac{\rho_\alpha}{(\alpha+\beta+\gamma)/2}$, in view of $\tan\theta = \frac{\rho}{x}$.

$\overline{|AM|} = \gamma + \overline{|BM|} = \gamma + \frac{\beta + \alpha - \gamma}{2} = \frac{\alpha + \beta + \gamma}{2}$. and from right triangle $A\Delta I$, we also have

Therefore we obtain, $\frac{\rho_\alpha}{(\alpha+\beta+\gamma)/2} = \frac{\rho}{x}$. Using the fact $x = \frac{\beta+\gamma-\alpha}{2}$ and $\rho = \frac{2(\text{AREA})}{\alpha+\beta+\gamma}$;

where AREA simply stands for the area of triangle $A\overset{\Delta}{B}\Gamma$ we finally arrive at

$2\rho_\alpha = \frac{4(\text{AREA})}{-\alpha+\beta+\gamma}$ by cyclicity or permutation of the three letters $A\overset{\Delta}{B}\Gamma$, we can state the

formulas for the four dameters:

$$2\rho = \frac{4(\text{AREA})}{\alpha+\beta+\gamma},\quad 2\rho_\alpha = \frac{4(\text{AREA})}{\beta+\gamma-\alpha},\quad 2\rho_\beta = \frac{4(\text{AREA})}{\alpha-\beta+\gamma},\quad 2\rho_\gamma = \frac{4(\text{AREA})}{\alpha+\beta-\gamma} \qquad (1)$$

Figure 1

$|\overline{AB}| = \gamma, |\overline{B\Gamma}| = \alpha, |\overline{\Gamma A}| = \beta$

$|\overline{A\Delta}| = x = |\overline{AE}|$

$|\overline{B\Delta}| = y = |\overline{BZ}|$

$|\overline{\Gamma Z}| = z = |\overline{\Gamma E}|$

$|\overline{\Delta I}| = |\overline{IZ}| = |\overline{IE}| = \rho$

$|\overline{O_\alpha M}| = \rho_\alpha$

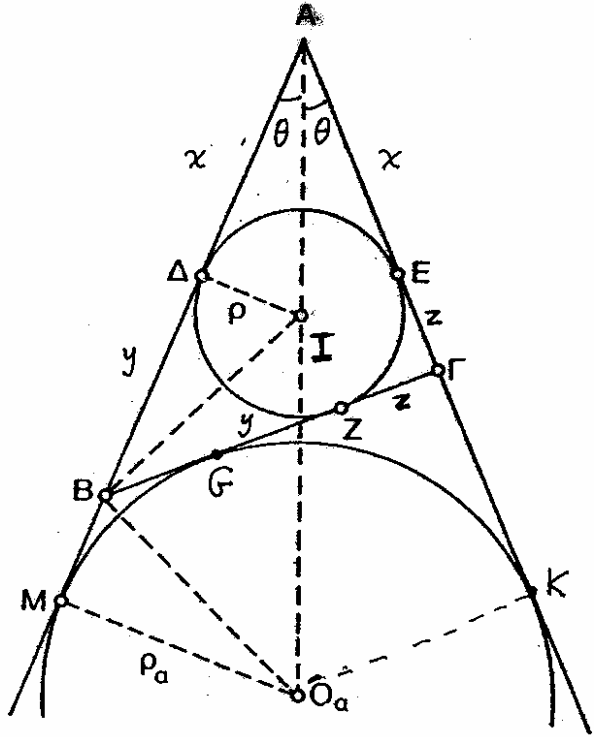

### 3. The Case of Right Triangles

If $\stackrel{\Delta}{AB\Gamma}$ is a right triangle with the ninety degree angle at the vertex A; then $\alpha$ becomes the hypotenuse length and $\beta$ and $\gamma$ the two leg lengths. We have,

$\alpha^2 = \beta^2 + \gamma^2 \Leftrightarrow 2\beta\gamma = (\beta+\gamma+\alpha)(\beta+\gamma-\alpha)$, which when combined with first two formulas in (1); and the fact that $4(\text{AREA}) = 2\beta\gamma$ ; it yields $2\rho = \beta + \gamma - \alpha$ and $2\rho_\alpha = \alpha + \beta + \gamma$

Note that $2\rho_\alpha$ is equal to the perimeter of the right triangle. Similarly, in view of the algebraic equivalence $2\beta\gamma = (\alpha+\beta-\gamma)[\alpha-(\beta-\gamma)]$; couched with (again) $4(\text{AREA}) = 2\beta\gamma$ , and the last two formulas in (1) ; we obtain, $2\rho_\beta = \alpha + \beta - \gamma$ and $2\rho_\gamma = \alpha + \gamma - \beta.$ Altogether,

$$2\rho = \beta + \gamma - \alpha,\ 2\rho_\alpha = \alpha + \beta + \gamma,\ 2\rho_\beta = \alpha + \beta - \gamma,\ 2\rho_\gamma = \alpha - \beta + \gamma \quad (2)$$

## 4. The Case of Pythagorean Triangles

The entire family of Pythagorean Triangles or triples $(\alpha, \beta, \gamma)$ can be described by the well-known parametric formulas,

$$\alpha = \delta(m^2 + n^2), \beta = \delta(2mn), \gamma = \delta(m^2 - n^2) \quad (3)$$

where m,n, $\delta$ are positive integers such that (m,n)=1, m>n, and $m + n \equiv 1 \pmod 2$ (i.e. m and n have different parities, one is even, the other odd). This parametric description can be found in most books of number theory; textbooks or more scholarly books. For example, you may refer to [1], [2], or [3].

For $\delta = 1$, the Pythagorean triple or triangle is called *primitive*; with hypotenuse length $\alpha = m^2 + n^2$, leg of even length $\beta = 2mn$, and leg of odd length $\gamma = m^2 - n^2$; with (m,n)=1, (4) $m + n \equiv 1$, $m, n \in \mathbb{Z}^+$, and $m > n$.

Applying formulas (2) and (4) we see that,

If $(\alpha, \beta, \gamma)$ is a primitive Pythagorean triple, the four diameters $2\rho, 2\rho_\alpha, 2\rho_\beta, 2\rho_\gamma$ are given by the formulas,

$$2\rho = 2n(m - n),\ 2\rho_\alpha = 2m(n + m),\ 2\rho_\beta = 2n(m + n),\ 2\rho_\gamma = 2m(m - n); \quad (5)$$

where the positive integers m,n satisfy the conditions in (4).

## 5. Seven Results from Number Theory

Result1: The entire set of solutions {x,y,z}, in positive integers, to the diophantine equation

$x^2 + 2y^2 = z^2$, with (x,y)=1, can be parametrically described by the formulas,

$x = |k^2 - 2\lambda^2|$, $y = 2\kappa\lambda$, $z = k^2 + 2\lambda^2$, where the parameters $k, \lambda$, can take any positive integer values which are relatively prime; $(k, \lambda) = 1, k, \lambda \in \mathbb{Z}^+$ and with k odd; $k \equiv 1 \pmod{2}$ (so that both x and z are odd) Result 1 is well-known in the literature and can be found in [1].

<u>Result 2:</u> The entire set of solutions in positive integers, to the diophantine equation

$x^2 + y^2 = 2z^2$, with (x,y)=1, can be parametrically (up to symmetry; the equation is symmetric with respect to x and y; if {a,b,c} is a solution, so is {b,a,c}) described by the formulas,

$x = |k^2 + 2k\lambda - \lambda^2|$, $y = |-k^2 + 2k\lambda + \lambda^2|$, $z = k^2 + \lambda^2$;

where $k, \lambda$ can be any positive integers such that $(k, \lambda) = 1$, $k + \lambda \equiv 1 \pmod{2}$ (so all three x,y,z, in any solution with (x,y)=1; must be odd). For historical commentary on the diophantine equation of result2 as well as a brief sketch of the derivation of the general solution, go to Section 11.

<u>Result 3:</u> Let $a, b \in \mathbb{Z}$, $ab \neq 0$, and (a,b)=1. The following statements hold true:

(i) If $c \in \mathbb{Z}$, and $a | b \cdot c$, then a must be a divisor of c.

(ii) If $p | c_1 c_2 \cdots c_k$; where p is a prime and $c_1, c_2, \ldots, c_k$ integers; then p must be a divisor of at least one of the numbers $c_1, c_2, \ldots, c_k$.

(iii) If $a \equiv b \equiv 1 \pmod{2}$, then $(a - b, a + b) = 2$ and either $a - b \equiv 0$ and $a + b \equiv 2 \pmod{4}$; or (vice-versa) $a - b \equiv 2$ and $a - b \equiv 0 \pmod{4}$

(iv) If $a + b \equiv 1 \pmod{2}$, then $(a^2 \pm 2ab \pm b^2, a^2 + b^2) = 1$ for any of the four combinations of signs, and if $a \equiv b \equiv 1 \pmod{2}$, then $(a^2 \pm 2ab \pm b^2, a^2 + b^2) = 2$, for any sign combination.

Parts (iii)-(v) of Result 3, are typically assigned as straightforward exercises in the first half of an introductory number theory course. Part (i) can be found in number theory books either listed as an exercise or as a theorem: see [2] or [3] . And part (ii) is used in proving the Fundamental

Theorem of Arithmetic (unique factorization of integers into prime powers). Again standard material found in number theory books.

Result 4: $a, b, c \varepsilon \mathbb{Z}^+, (a,b) = 1$, and $ab = c^n$; then $a = c_1^n, b = c_2^n$, for some $c_1, c_2 \varepsilon \mathbb{Z}^+$, with $(c_1, c_2) = 1$, and $c_1 c_2 = c$

This is a well-known result, you may refer to [2] or [3]. Based on Result 4, the following result can be easily proven.

Result 5: If p is a prime number, and $ab = pc^n$; $a, b, c \varepsilon \mathbb{Z}^+, (a,b) = 1$; then there exist $c_1, c_2 \varepsilon \mathbb{Z}^+$ such that either $a = pc_1^n, b = c_2^n$ or alternatively $a = c_1^n$ and $b = pc_2^n$; and with $(c_1, c_2) = 1$ and $c_1 c_2 = c$

The proof of Result 5 is easy, if we make use of Results 4 and 3(ii). Indeed, by Result 3(ii), we know that p must divide at least one of a and b; say $p|a$. We then have a=pd, for some $d \varepsilon \mathbb{Z}^+$. Accordingly, $ab = pc^n$ gives $pdb = pc^n$; $db = c^n$. Clearly (d,b)=1, since d is a divisor of a and (a,b)=1. By Result 4 we must have $d = c_1^n$, $b = c_2^n$, with $(c_1, c_2) = 1$ and $c_1 c_2 = c$. Thus, $a = pc_1^n$ and $b = c_2^n$. The other possibility, $p|b$, leads to $a = c_1^n$ and $b = pc_2^n$.

Result 6: (i) If p is a prime, and $a, b, c \in \mathbb{Z}^+$, with $(a,b) = 1$ and $pab = c^n$; then either $a = p^{n-1} \bullet c_1^n$ and $b = c_2^n$; or $a = c_1^n$ and $b = p^{n-1} \bullet c_2^n$; for some positive integers $c_1, c_2$ such that $(c_1, c_2) = 1$ and $c_1 c_2 = c$

(ii) In particular when n=2 in part (i), we must have either $a = pc_1^2$ and $b = c_2^2$; or, alternatively $a = c_1^2$ and $b = pc_2^2$.

We leave the proof of Result 6 to the reader as an exercise. Here is a sketch of this proof. Apply Result 3 (ii) to conclude that the prime p must divide c. Then use this fact to deduce that $p^{n-1}$ must divide either a or b; by using the hypothesis(a,b)=1 together with the fact that $p^{n-1}$ is a prime power; and in conjunction with Result 3 (i). Finally, apply Result 4 to conclude the proof.

Result 7: If $a, b \in \mathbb{Z}^+$ and $a^n | b^n$, then $a | b$

In words, this can be stated as follows: If the nth power of a positive integer divides the nth power of another positive integer; then the former positive integer must divide the latter. This is a well-known result, refer to [2] or [3].

**6. Infinitely many primitive Pythagorean triangles with one leg length being a square.**

Let us remark at the outset, that a Pythagorean Triangle can have only one leg whose length is a square. P. Fermat was the first historically known mathematician who's used the method of infinite descent. He used this method (see Section 8) to show that the diophantine equation $x^4 + y^4 = z^2$ has no positive integer solutions. This is a well-known result which can be found in many number theory books. An immediate corollary of this is the non-existence of Pythagorean Triangles both of whose leg lengths are integer squares. As we see below, the process of describing all primitive Pythagorean Triangles with one leg length being a square, is a fairly simple one.

*Proposition 1*

*(i) There exist infinitely many primitive Pythagorean triangles or triples { $\alpha, \beta, \gamma$ }, with the even leg length $\beta$ being a square. The entire family can be parametrically described by the formulas,*

$\alpha = m^2 + n^2$, $\beta = 2mn$, $\gamma = m^2 - n^2$;, *where* $m = 2t_1^2$ *and* $n = t_2^2$; *or* $m = t_1^2$, *and* $n = 2t_2^2$, *where* $t_1$ *and* $t_2$ *can be any positive integers such that* $(t_1, t_2) = 1$ *and* $t_2$ *odd in the first case*; $t_1$ *odd in the second, and also with* $2t_1^2 > t_2^2$ *or* $t_1^2 > 2t_2^2$ *respectively.*

*(ii) There exist infinitely many primitive Pythagorean triangles or triples* $\{\alpha, \beta, \gamma\}$, *with the odd length* $\gamma$ *being a square. The entire family can be parametrically described by the formulas,*

$\alpha = m^2 + n^2$, $\beta = 2mn$, $\gamma = m^2 - n^2$; *where* $m = t_1^2 + t_2^2$, $n = 2t_1 t_2$, *where* $t_1$ *and* $t_2$ *can be any positive integers, such that* $(t_1, t_2) = 1$, $t_1 + t_2 \equiv 1 \pmod{2}$ ( $t_1$ *and* $t_2$ *have different parities*), *and* $t_1 > t_2$ *(this ensures that there is no repetition of triples* $\{\alpha, \beta, \gamma\}$ *)*

Proof: (i) If $\alpha, \beta, \gamma$ are given by the stated formulas, then an easy calculation shows that $\beta$ is square and that $\{\alpha, \beta, \gamma\}$ is a primitive Pythagorean triple. Now, conversely, suppose that $\{\alpha, \beta, \gamma\}$ is a primitive Pythagorean triple with $\beta$ a square. We must have $\alpha = m^2 + n^2$, $\beta = 2mn$, $\gamma = m^2 - n^2$, and $\beta = c^2$; for some positive integers m, n, c, such that (m, n)=1, m+n ≡ 1 (mod2), and m>n. From $2mn = c^2$ and Result 6 with p=2, it follows that either $m = 2t_1^2$ and $n = t_2^2$ or $m = t_1^2$ and $n = 2t_2^2$.

In the first case we have $m = 2t_1^2$, and $n = t_2^2$. Now, (m,n)=1 and m+n ≡ 1(mod2); which easily implies that $t_2$ is odd and $(t_1, t_2) = 1$. And m>n implies $2t_1^2 > t_2^2$. In the second case, $m = t_1^2$ and $n = 2t_2^2$. Again by virtue of (m,n)=1, m>n, and m+n ≡ 1(mod2), it easily follows that $t_1$ must be odd, $(t_1, t_2) = 1$, and $t_1^2 > 2t_2^2$.

(ii) As with part (i), one direction in the proof is almost immediate. In the other direction, suppose that $\{\alpha, \beta, \gamma\}$ is a primitive Pythagorean triple with the odd length $\gamma$ being a square. We must have, $\alpha = m^2 + n^2, \beta = 2mn, \gamma = m^2 - n^2$, and $\gamma = c^2$; for some positive integers m,n such that (m,n=1, m+n ≡ 1 (mod2) and m>n. Then, $\gamma = c^2 \Leftrightarrow m^2 = n^2 + c^2$. Clearly (m,n)=1, easily implies from the last equation that also (m,c)=1=(n,c). In other words {m,n,c} is a primitive Pythagorean triple with m being the hypotenuse length and hence odd; and also from m+n ≡ 1(mod2), we infer that n must be even. Consequently, $m = t_1^2 + t_2^2, n = 2t_1 t_2, c = t_1^2 - t_2^2$, for some $t_1, t_2 \varepsilon \mathbb{Z}_-^+$, such that $(t_1, t_2) = 1, t_1 + t_2 \equiv 1 (\mod 2), \text{and } t_1 > t_2$ $\square$

**7. Primitive Pythagorean triangles with one leg length being a square and one of the diameters $2\rho, 2\rho_\alpha, 2\rho_\beta, 2\rho_\gamma$ also a square.**

If $\{\alpha, \beta, \gamma\}$ is a primitive Pythagorean triple with

$\alpha = m^2 + n^2, \beta = 2mn, \gamma = m^2 - n^2, m, n \varepsilon \mathbb{Z}^+, (m, n) = 1, m + n \equiv 1 (\mod 2)$, and m>n, then there are exactly eight combinations (with one side length a square and one of the four diameters also a square.)

Combination 1: $\beta$ a square and $2\rho$ a square

Combination 2: $\beta$ a square and $2\rho_\alpha$ a square

Combination 3: $\beta$ a square and $2\rho_\beta$ a square

Combination 4: $\beta$ a square and $2\rho_\gamma$ a square

Combination 5: $\gamma$ a square and $2\rho$ a square

Combination 6: $\gamma$ a square and $2\rho_\alpha$ a square

Combination 7: $\gamma$ a square and $2\rho_\beta$ a square

Combination 8: $\gamma$ a square and $2\rho_\gamma$ a square.

As Theorem 1 below shows, the two combinations 6 and 8, cannot really occur. On the other hand, according to Theorem 2, the remaining combinations do occur; and in each case (combination) the entire family of such triples is parametrically described.

**Theorem 1**

(i) There are no primitive Pythagorean Triangles with the odd length $\gamma$ a square and with diameter $2\rho_\alpha$ (also the perimeter of the triangle) also a square.

(ii) There are no primitive Pythagorean Triangles with the odd length $\gamma$ a square, and the diameter $2\rho_\gamma$ also a square.

**Proof**

(i) By Proposition 1(ii), if such a triple $\{\alpha, \beta, \gamma\}$ exists, we must have

$\alpha = m^2 + n^2, \beta = 2mn, \gamma = m^2 - n^2$, and with $m = t_1^2 + t_2^2, n = 2t_1 t_2, (t_1, t_2) = 1, t_1 + t_2 \equiv 1 \pmod{2}$, for some $t_1, t_2 \varepsilon \mathbb{Z}^+$

In addition, $2\rho_\alpha = K^2$, for some positive integer K. By (5), we know that $2\rho_\alpha = 2m(m+n)$. Therefore, $2(t_1^2 + t_2^2)(t_1^2 + t_2^2 + 2t_1 t_2) = K^2$; which is contradictory since, in view of $t_1 + t_2 \equiv 1 \pmod{2}$, the left-hand side is congruent to 2; but the right-hand side can only be congruent to 0 or 1 mod 4.

(ii) As in (ii), we must have $\alpha = m^2 + n^2, \beta = 2mn, \gamma = m^2 - n^2, m = t_1^2 + t_2^2, n = 2t_1 t_2$; for some positive integers $t_1, t_2$ such that $(t_1, t_2) = 1$ and $t_1 + t_2 \equiv 1 \pmod{2}$. Additionally, $2\rho_\gamma = N^2$, for some $N \varepsilon \mathbb{Z}^+$. And from (5), $2\rho_\gamma = 2m(m-n)$. Consequently, $2(t_1^2+t_2^2)(t_1^2+t_2^2-2t_1t_2) = N^2$;

$2(t_1^2+t_2^2)(t_1-t_2)^2 = N^2$, which is impossible modulo 4, since the left-hand side is congruent to 2, while the right-hand side can only be 0 or 1 mod 4.

### Theorem 2

(i) The primitive Pythagorean triples $\{\alpha, \beta, \gamma\}$ with $\beta$ (the even length) a square and the diameter $2\rho_\alpha$ (also the perimeter) also a square, are precisely the members of the following family:

$$\text{Family } F_1: \left\{ \begin{array}{l} \alpha = m^2+n^2, \beta = 2mn, \gamma = m^2-n^2, \text{and with} \\ m = 2t^2, n = t_2^2, 2t_1^2 > t_2^2, \text{where} \\ t_2 = |\kappa^2 - 2\lambda^2|, t_1 = 2\kappa\lambda, \kappa \equiv 1 (\text{mod } 2), \\ (\kappa, \lambda) = 1, \kappa, \lambda \varepsilon \mathbb{Z}^+ \end{array} \right\}$$

(ii) The primitive Pythagorean triples $\{\alpha, \beta, \gamma\}$ with $\beta$ (the even length) a square and the diameter $2\rho_\beta$ also a square, are precisely the members of the following family:

$$\text{Family } F_2: \left\{ \begin{array}{l} \alpha = m^2+n^2, \beta = 2mn, \gamma = m^{-2}n^2, \text{and with} \\ m = t_1^2, n = 2t_2^2, t_1^2 > 2t_2^2, \text{where} \\ t_1 = |\kappa^2 - 2\lambda^2|, t_2 = 2\kappa\lambda, \kappa \equiv 1 (\text{mod } 2), \\ (\kappa, \lambda) = 1, \kappa, \lambda \varepsilon \mathbb{Z}^+ \end{array} \right\}$$

(iii) The primitive Pythagorean triples $\{\alpha, \beta, \gamma\}$ with $\beta$ (the even length) a square and the diameter $2\rho_\gamma$ also a square, are precisely the members of the following family:

$$\text{Family } F_3: \left\{ \begin{array}{l} \alpha = m^2+n^2, \beta = 2mn, \gamma = m^2-n^2, \text{with} \\ m = 2t_1^2, n = t_2^2, 2t_1^2 > t_2^2; \text{and with} \\ t_1 = \kappa^2 + \lambda^2, t_2 = |-\kappa^2 + 2k\lambda + \lambda^2| \\ \kappa, \lambda \varepsilon \mathbb{Z}^+, \kappa + \lambda \equiv 1 (\text{mod } 2), (\kappa, \lambda) = 1 \end{array} \right\}$$

(iv) The primitive Pythagorean triples $\{\alpha, \beta, \gamma\}$ with $\gamma$ (the odd length) a square and the inner diameter $2\rho$ (also the perimeter) also a square are precisely the members of the following family:

Family $F_4$: $\left\{ \begin{array}{l} \alpha = m^2 + n^2, \beta = 2mn, \gamma = m^2 - n^2, \text{and with} \\ m = t_1^2 + t_2^2, n = 2t_1 t_2; \text{where} \\ t_1 = \kappa^2, t_2 = \lambda^2, t_1 > t_2 \\ \kappa, \lambda \varepsilon \mathbb{Z}^+, (\kappa, \lambda) = 1, \kappa + \lambda \equiv 1 (\text{mod } 2) \end{array} \right\}$

(v) The primitive Pythagorean triples $\{\alpha, \beta, \gamma\}$ with $\gamma$ (the odd length) a square and the diameter $2\rho_\beta$ also a square, are precisely the members of the family $F_4$ of part (iv) of this Theorem. This is because (as it will be easily seen in the proof) when $\gamma$ is square; then $2\rho$ is also a square in, and only if, $2\rho_\beta$ is also a square.

(vi) The primitive Pythagorean triples $\{\alpha, \beta, \gamma\}$ with $\beta$ (the even length) a square and the inner diameter $2_\rho$ also a square, are precisely the members of the following family:

Family $F_6$: $\left\{ \begin{array}{l} \alpha = m^2 + n^2, \beta = 2mn, \gamma = m^2 - n^2, \text{and with} \\ m = t_1^2, n = 2t_2^2; \text{where,} \\ t_1 = \kappa^2 + 2\lambda^2, t_2 = 2\kappa\lambda; \kappa, \lambda \varepsilon \mathbb{Z}^+, \text{and} \\ \text{such that } \kappa \equiv 1 (\text{mod } 2) \text{ and } (\kappa, \lambda) = 1 \\ \text{and with } t_1^2 > 2t_2^2 \end{array} \right\}$

**Proof**

First, we point out that in each part, one direction in the proof is straightforward: if $\{\alpha, \beta, \gamma\}$ is a member of the given family, then it is a primitive Pythagorean triple with the stated property. Below, in each part, we prove the converse. If $\{\alpha, \beta, \gamma\}$ is a primitive Pythagorean triple with the stated property, then it must be a member of the given family.

(i) By Proposition 1(i) we know that we must have, $\alpha = m^2 + n^2$, $\beta = 2mn$; and with either $m = 2t_1^2$, $n = t_2^2$, and $t_2 \equiv 1$; or alternatively $m = t_1^2$, $n = 2t_2^2$, and $t_1 \equiv 1 \pmod 2$. And in either case, with $(t_1, t_2) = 1$. By (5) we know that, $2m(m+n) = L^2$, for some $L \in \mathbb{Z}^+$. Since m and n have different parities, $m + n \equiv 1 \pmod 2$. But obviously, L must be even; and so $L^2 \equiv 0 \pmod 4$. Thus, $2m(m+n) \equiv 0 \pmod 4$; which implies that m is even, since m+n is odd. Therefore, only the first of the above two possibilities can occur. Namely, $m = 2t_1^2$ and $n = t_2^2$, $(t_1, t_2) = 1$ and with $t_2$ odd (and also with $2t_1^2 > t_2^2$) We have

$2m(m+n) = L^2 \Leftrightarrow 4t_1^2 \cdot (2t_1^2 + t_2^2) = L^2 \Rightarrow (2t_1)^2 | L^2 \Rightarrow 2t_1 | L$, by Result 7. We set $L = (2t_1) \cdot L_1$, for some positive integer $L_1$. Furthermore, $4t_1^2(2t_1^2 + t_2^2) = 4t_1^2 L_1^2 \Rightarrow t_2^2 + 2t_1^2 = L_1^2$; which shows that $\{t_2, t_1, L_1\}$ is a solution in positive integers to the diophantine equation $x^2 + 2y^2 = z^2$. By Result 1 we must have, $t_2 = |k^2 - 2\lambda^2|$, $t_1 = 2\kappa\lambda$, $\kappa \equiv 1 \pmod 2$, $(\kappa, \lambda) = 1$, for some positive integers $\kappa, \lambda$.

(ii) This derivation is very similar to the one in part (i). Indeed, since by (5) we must have $2\rho_\beta = 2n(m+n) = L^2$, for some $L \in \mathbb{Z}^+$; a congruence modulo 4 implies that since m+n is odd, n must be even. This couched with Proposition 1(i), yields $\alpha = m^2 + n^2, \beta = 2mn, m = t_1^2, n = 2t_2^2$, with $t_1$ odd, $(t_1, t_2) = 1$, and $t_1^2 > 2t_2^2$. Then, from $2n(m+n) = L^2$ we easily obtain

$(2t_2)^2 \cdot (t_1^2 + 2t_2^2) = L^2$ and by applying Result 7, we end up with $t_1^2 + 2t_2^2 = L_1^2$, for some $L_1 \in \mathbb{Z}^+$;

which shows that $\{t_1, t_2, L_1\}$ is a solution in positive integers to the diophantine equation

$x^2 + 2y^2 = z^2$; hence, $t_1 = |\kappa^2 - 2\lambda^2|, t_2 = 2\kappa\lambda$, for some positive integers $\kappa, \lambda$ such that

$(\kappa, \lambda) = 1$, and $\kappa \equiv 1 \pmod{2}$.

(iii)   Since $2\rho_\gamma = 2m(m-n)$ is a square and m-n is odd (by virtue of $m+n \equiv 1 \pmod{2}$), it

follows that m must be even. (Again arguing modulo 4). Combining this with Proposition 1 (i)

gives $m = 2t_1^2, n = t_2^2$, with $t_2$ odd, $2t_1^2 > t_2^2$, and $(t_1, t_2) = 1$. Then, from $\rho_\gamma$=square we obtain,

$(2t_2)^2 \cdot (2t_1^2 - t_2^2) = L^2$, for some positive integer. By Result 7 it follows that $2t_2$ is a divisor of L;

put $L = (2t_2) \cdot L_1$, for some integer $L_1$. Therefore the last equation gives,

$2t_1^2 - t_2^2 = L_1^2; 2t_1^2 = t_2^2 + L_1^2$. Note that since $t_2$ is odd, so must be $L_1$; and $t_1$ must be odd as well,

as a congruence mod4 easily shows. Therefore $\{t_2, L_1, t_1\}$ is a positive integer solution to the

diophantine equation $x^2 + y^2 = 2z^2$, with (x,y)=1. Note that $(t_1, t_2) = 1$ easily implies

$(t_2, L_1) = 1 = (t_1, t_2)$ in the above equation. Accordingly by Result 2 we must have either

$t_1 = \kappa^2 + \lambda^2$, $t_2 = |-\kappa^2 + 2\kappa\lambda + \lambda^2|$; or alternatively $t_1 = \kappa^2 + \lambda^2, t_2 = |-\kappa^2 - 2\kappa\lambda + \lambda^2|$; or positive

integers $\kappa, \lambda$ such that $(\kappa, \lambda) = 1$ and $\kappa + \lambda \equiv 1 \pmod{4}$.

(iv)   Since $2\rho$ must be a square, by (5) we know that $2\rho = 2n(m-n) =$ square; and since m-n

is odd, we must have $2n \equiv 0 \pmod{4}$; which means that n must be even.

Additionally, by Proposition 1(ii), we must have $m = t_1^2 + t_2^2$ and $n = 2t_1t_2$, with

$(t_1, t_2) = 1 t_1 + t_2 \equiv 1 \pmod{2}$, and $t_1 > t_2$. Altogether, 2n(m-n)=square

$\Leftrightarrow 4t_1t_2(t_1^2 + t_2^2 - 2t_1t_2) = L^2 \Leftrightarrow t_1t_2 \cdot [2(t_1 - t_2)]^2 = L^2$; which implies that, by Result 7, the positive integer $2(t_1 - t_2)$ is a divisor of L; and so, $L = 2(t_1 - t_2) \cdot L_1$; going back, $t_1 t_2 \cdot [2(t_1 - t_2)]^2 = 4(t_1 - t_2)^2 \cdot L_1^2$; $t_1 t_2 = L_1^2$, and since $(t_1, t_2) = 1$, Result 4 (applied with n=2) tells us that $t_1 = \kappa^2$ and $t_2 = \lambda^2$, for some positive integers $\kappa, \lambda$ such that $(\kappa, \lambda) = 1$. But, by virtue of $t_1 + t_2 \equiv 1 \pmod 2$, we must also have $\kappa + \lambda \equiv 1 \pmod 2$.

(v) Here, we only have to observe that since $\gamma$ is a square, we must have, by Proposition 1(ii), $\alpha = m^2 + n^2, \beta = 2mn, \gamma = m^2 - n^2$, $m = t_1^2 + t_2^2, n = 2t_1 t_2, (t_1, t_2) = 1, t_1 + t_2 \equiv 1 \pmod 2$, and $2\rho_\beta = 2n(m+n) = 4t_1 t_2 \cdot (t_1 + t_2)^2$ and 2n (m-n)$= 4t_1 t_2 (t_1 - t_2)^2 = 2\rho$. And since $(t_1, t_2) = 1$, it is then clear that, $2\rho_\beta =$ square $\Leftrightarrow 2\rho =$ square $\Leftrightarrow (t_1 = \kappa^2 \text{and} t_1 = \lambda^2)$, for some positive integers $\kappa, \lambda$, with $(\kappa, \lambda) = 1 \text{and} \kappa + \lambda \equiv 1 \pmod 2$.

(vi) From formulas (5), we must have $2\rho = 2n(m - n) =$ square; which as we have seen previously in other parts of this Theorem's proof, easily implies that since m-n is an odd integer; n must be even. But then $\beta$ being a square implies (by Proposition 1(i)) that, $m = t_1^2, n = 1 t_2^2, t_1 \equiv 1 \pmod 2, (t_1, t_2) = 1, \text{and} t_1^2 > 2 t_2^2$. We have, 2n(m-n)$= L^2$, for some $L \in \mathbb{Z}^+$; thus, $(2t_2)^2 \cdot (t_1^2 - 2t_2^2) = L^2$; which implies (by Result 7) that $2t_2$ is a divisor of L; and so $L = (2t_2) \cdot L_1$, for some positive integer $L_1$. Going back, $(2t_2)^2 \cdot (t_1^2 - 2t_2^2) = (2t_2)^2 \cdot L_1^2 \Leftrightarrow L_1^2 + 2t_2^2 = t_1^2$, which shows that $\{L_1, t_2, t_1, \}$ is a positive integer solution to the diophantine equation $x^2 + 2y^2 = z^2$, with (x,y)=1 (note that $(t_1, t_2) = 1 \Rightarrow (L_1, t_2) = 1$) Hence, according to Result 1 we must have, $t_1 = \kappa^2 + 2\lambda^2, t_2 = 2\kappa\lambda$, for some positive integers $\kappa \text{and} \lambda$, with $\kappa$ odd and $(\kappa, \lambda) = 1$. This concludes the proof..

## Numerical Examples

1) For Family $F_1$.

<u>Example 1</u>: Take

$\kappa = 1, \lambda = 1;$ then $t_2 = 1, t_1 = 2, m = 8, n = 1, \alpha = 65, \beta = (4)^2 = 16, \gamma = 63, 2\rho_\alpha = (12)^2 = 144$

<u>Example 2</u>: Take

$\kappa = 3, \lambda = 2;$ then $t_2 = 1, t_1 = 12, m = 288, n = 1, \alpha = 82945, \beta = (24)^2 = 576, \gamma = 82943, 2\rho_\alpha = (408)^2 = 166464$

2) For Family $F_2$

<u>Example 3</u>: Take

$\kappa = 1, \lambda = 2,$ then $t_1 = 7, t_2 = 4, m = 49, n = 32, \alpha = 3425, \beta = (58)^2 = 3364, \gamma = 1377, 2\rho_\beta = (144)^2 = 20736$

<u>Example 4</u>: Take

$\kappa = 5, \lambda = 1;$ then $t_1 = 23, t_2 = 10, m = 529, n = 200, \alpha = 319841, \beta = (460)^2 = 211600, y = 239841, 2\rho_\beta = (540)^2 = 291600$

3) For Family $F_3$

<u>Example 5</u>: Take

$\kappa = 1, \lambda = 2;$ then $t_1 = 5, t_2 = 7, m = 50, n = 49, \alpha = 4901, \beta = (70)^2 = 4900, \gamma = 99, 2\rho_\gamma = (10)^2 = 100$

<u>Example 6</u>: Take

$\kappa = 2, \lambda = 1;$ then $t_1 = 5, t_2 = 1, m = 50, n = 1, \alpha = 2501, \gamma = 2499, \beta = (10)^2 = 100, 2\rho_\gamma = (70)^2 = 4900$

<u>Example 7</u>: Take $\kappa = 2, \lambda = 3$ (use the alternative formula for $t_2$); then

$t_1 = 13, t_2 = 17, m = 338, n = 289, \alpha = 197765, \beta = (442)^2 = 195364, \gamma = 30723, 2\rho_\gamma = (182)^2 = 33124$

4) For Family $F_4 = F_5$

Example 9: Take

$\kappa = 1, \lambda = 2;$ then $t_1 = 1, t_2 = 4, m = 17, n = 8, \alpha = 353, \beta = 272, \gamma = (15)^2 = 225, 2\rho = (12)^2 = 144, 2\rho_\beta = (20)^2 = 400$

Example 10: Take

$\kappa = 2, \lambda = 3;$ then $t_1 = 4, t_2 = 9, m = 97, n = 72, \alpha = 14593, \beta = 13968, \gamma = (65)^2 = 4225, 2\rho = (60)^2 = 3600, 2\rho_\beta = (156)^2 = 24336$

5) For Family $F_6$

Example 11: Take

$\kappa = 1, \lambda = 1;$ then $t_1 = 3, t_2 = 2, m = 9, n = 8, \alpha = 97, \beta = (12)^2 = 144, \gamma = 17, 2\rho = (4)^2 = 16$

Example 12: Take

$\kappa = 1, \lambda = 2;$ then $t_1 = 9, t_2 = 4, m = 81, n = 32, \alpha = 7585, \beta = (72)^2 = 5384, \gamma = 5537, 2\rho = (56)^2 = 3136$

8. **The diophantine equation $x^2+y^2=2z^2$ and closing remarks.**

A. *A brief description / sketch of the derivation of the general solution of the equation*

$x^2+y^2=2z^2$ *with* $(x,y)=1$

Below we outline a method for obtaining all the positive integer solutions to this equation, via two-parameter formulas. To this end, let {a,b,c} be such a solution. We must then have,

$$\begin{cases} a^2+b^2=2c^2 \\ (a,b)=1; \text{ and } a,b,c \in \mathbb{Z}^+ \end{cases}$$

Obviously, because of the condition (a,b)=1; either a and b are both odd or otherwise they must have different parities. But $2c^2 \equiv 0 \text{ or } 2$ modulo 4 (for c even or odd respectively). If a and b had different parity then, $a^2+b^2 \equiv 1 \pmod{4}$. This clearly shows that only the first possibility can occur, namely a and b being both odd; and so c being odd as well. Next, observe that the condition (a,b)=1 easily implies (b,c)=1=(a,c); from which it easily follows that any two of the three integers a,b,c, must be distinct; with a single exception, the case a=b=c=1. The above equation is the equivalent to $(\frac{a}{c})^2+(\frac{b}{c})^2=2$, which shows that the point $(a/c, b/c)$ on the two-dimensional XY plane lies on the circle of center the origin (0,0) and radius $\sqrt{2}$:

The circle with equation $X^2+Y^2=2$.

The point $(a/c, b/c)$ is a rational point (i.e. a point both of whose coordinates are rational numbers) on the circle. Assume that at least one of the positive integers a, b, c is greater than 1. Which means (according to the above explanation) that all three of them are pairwise distinct and greater than 1. Then the point (1,1) lies on the same circle and it is distinct from $(a/c, b/c)$.

Consider the straight line on the XY plane that passes through these two points. It must have slope $s = \dfrac{b/c - 1}{a/c - 1}$; with $s \neq 0$. Since the point $(a/c, b/c)$ lies on the above circle and in the first quadrant of the XY plane; it follows that $s<0$.

The slope s is a negative rational number; so that s can be represented as $s = -\dfrac{K}{\lambda}$, where K and $\lambda$ are positive integers, with $(K, \lambda) = 1$. We have, $s = -\dfrac{K}{\lambda} = \dfrac{b/c - 1}{a/c - 1}$.

If we solve for $b/c$ in terms of $a/c$, and $K/\lambda$ and substitute in the circle equation

$(a/c)^2 + (b/c)^2 = 2$; after a certain amount of algebra we arrive at,

$\left[(\dfrac{K}{\lambda})^2 + 1\right] \cdot (a/c)^2 - 2(\dfrac{K}{\lambda}) \cdot (\dfrac{K}{\lambda} - 1)(\dfrac{a}{c}) + (\dfrac{K}{\lambda})^2 - 2(\dfrac{K}{\lambda}) - 1 = 0.$

The last equation shows that the number $a/c$ is one of the two roots of the quadratic equation (in the variable t), $\left[(\dfrac{K}{\lambda})^2 + 1\right] \cdot t^2 - 2\left(\dfrac{K}{\lambda}\right)\left(\dfrac{K}{\lambda} - 1\right) \cdot t + (\dfrac{K}{\lambda})^2 - 2(\dfrac{K}{\lambda}) - 1 = 0.$

The other root is the number 1, as it can be seen by inspection. But the sum of the two roots of a quadratic equation $\alpha t^2 + \beta t + \gamma = 0$, must equal $-\beta/\alpha$. Applying this fact to the above quadratic equation whose roots are 1 and $a/c$; and solving for $a/c$ we obtain, $\dfrac{a}{c} = \dfrac{K^2 + 2K\lambda - \lambda^2}{K^2 + \lambda^2}$

Moreover, going back to the equation for the slope s and solving for $b/c$ we also obtain,

$\dfrac{b}{c} = \dfrac{-K^2 + 2K\lambda + \lambda^2}{K^2 + \lambda^2}$

At this point, the following simple result from number theory may be used: If $a_1, a_2, a_3, a_4$, are four positive integers such that $\frac{a_1}{a_2} = \frac{a_3}{a_4}$ and $(a_1, a_2) = 1$; then $a_3 = \delta a_1$ and $a_4 = \delta a_2$; where $\delta = (a_3, a_4)$ (The greatest common divisor of $a_3$ and $a_4$). This typically can be assigned as a simple exercise in elementary number theory, which can be easily proved with the aid of Result 3(i). In our situation we have $\frac{a}{c} = \frac{K^2 + 2K\lambda - \lambda^2}{K^2 + \lambda^2}$ and $\frac{b}{c} = \frac{-K^2 + 2K\lambda + \lambda^2}{K^2 + \lambda^2}$

Note that in fact, since a,b,c, and $K^2 + \lambda^2$ are all positive; both integers $K^2 + 2K\lambda - \lambda^2$ and $-K^2 + 2K\lambda + \lambda^2$ must also be positive. Also recall that $(K, \lambda) = 1$

Now, by Result 3(iv), if $\delta = (K^2 + 2K\lambda - \lambda^2, K^2 + \lambda^2)$, then $\delta = 1$ or $2$, respectively for $K + \lambda \equiv 1 \pmod{2}$; or $K \equiv \lambda \equiv 1 \pmod{2}$. But if $\delta = 2$, then according to the above fact both integers a and c would have to be multiples of 2, which obviously cannot happen, since a and c are both odd. The same argument applies to the second fraction; which leads to the conclution that, $a = K^2 + 2K\lambda - \lambda^2, b = -K^2 + 2K\lambda + \lambda^2, c = K^2 + \lambda^2$; and with $(K, \lambda) = 1$ and $K + \lambda \equiv 1 \pmod{2}$. Conversely, if a,b,c satisfy these parametric formulas, then a routine calculation shows that $a^2 + b^2 = 2c^2$. Even though some dublication of solutions may occur; to ensure that for any choice of positive integers K and $\lambda$, the above formulas produce positive a and b (c is always positive), we may take $a = |K^2 + 2K\lambda - \lambda^2|, b = |-K^2 + 2K\lambda + \lambda^2|, c = K^2 + \lambda^2$.

In conclusion, the entire family of positive integer solutions to the diophantine equation

$\begin{cases} x^2 + y^2 = 2z^2 \\ (x, y) = 1 \end{cases}$ can be described by,

$x = \left| K^2 + 2K\lambda - \lambda^2 \right|, y = \left| -K^2 + 2K\lambda + \lambda^2 \right|, z = K^2 + \lambda^2$, where K, $\lambda$ can be any positive integers;

such that $(K, \lambda) = 1$, $K + \lambda \equiv 1 \pmod{2}$.

Note that the solution {1,1,1} is the only one which cannot be obtained from the above formulas.

*B. Historical Remarks*

According to L.E. Dickson's book (see[1]), historically, the diophantine equation $x^2 + y^2 = 2z^2$ has been investigated as a result of the problem of finding three integers whose squares are in arithmetic progression. A number of investigators are mentioned. Various researchers used a number of different techniques and methods to tackle this problem; most of them achieved partial solutions or families of solutions, but not the general solution (i.e. entire family of solutions). Below we mention five researchers who worked on this problem, as cited in Dickson's book. Another twenty or so names of individuals are mentioned as contributors to the understanding and solving of the above diophantine equation.

Diophantus (circa 150AD-250AD) is reported to have discovered three particular or special numbers whose squares are in arithmetic progression. In the thirteenth century, Jordanus Nemorarius found a particular family of solutions:

$x = b^2 - c^2/2$, $y = b^2 + 2bc + c^2/2$, $z = b^2 + bc + c^2/2$, where b can be any integer, while c is an even integer.

P. Fermat, in the 17$^{th}$ century, discovered the following family of solutions;

$x = r^2 - 2s^2$, $y = r^2 + 4rs + 2s^2$, $z = r^2 + 2rs + 2s^2$.

Then there is the 1918 paper of A.E. Jones, in which he discovered a family of integer triples, such that in each triple, the squares of the three numbers are in arithmetic progression; and the sum of any two of the three numbers in the triple is a square as well.

Finally, A. Desboves in the late 19$^{th}$ century, gave all the solutions to the more general diophantine equation, $x^2 + y^2 = (a^2 + b^2)z^2$ (in our case a=b=1)

**References**

**[1]** Dickson, L.E., *History of Theory of Numbers, Vol. II*, AMS Chelsea Publishing, Providence, Rhode Island, 1992 ISBN: )-8218-1935-6; 803 p.p. (unaltered textual reprint of original book, first published by Carnegie Institute of Washington in 1919, 1920, and 1923)

(a) For the diophantine equation $x^2 + 2y^2 = z^2$, pages; 420, 421, 426

(b) For the diophantine equation $x^2 + y^2 = 2z^2$, pages; 427, 428; 435-438

(c.) For Pythagorean triangles, pages 165-169; for more general information on rational right triangles, pages 170-190.

**[2]** Rosen, Kenneth H., *Elementary Number Theory and its Applications,* Third Edition, 1993, Addison-Wesley Publishing Company (there is now a fourth edition as well), 544 p.p

ISBN: 0-201-57889-1

(a) For Pythagorean triples, pages 436-442

(b) For Result 3 (i), page 91, Lemma 2.3

(c) For Result 4, page 101, exercise 57

(d) For Result 7, page 100, exercise 35

**[3]** Sierpinski, W., *Elementary Theory of Numbers,* originals edition, Warsaw, Poland, 1964, 480 p.p. (no ISBN numbers). More recent version (1988) published by Elsevier Publishing, and distributed by North-Holland. North-Holland Mathematical Library **32**, Amsterdam (1988) This book is available by various libraries, but it is only printed upon demand. Specifically, UMI Books on Demand, from: Pro Quest Company, 300 North Zeep Road, Ann Arbor, Michigan, 4801-1356 USA; ISBN: 0-598-52758-3

(a) For Result 3(i): pages 14,15

(b) For Result 4: pages 17,18

(c) For Result 7: page 15

**[4]** K. Zelator, *The diophantine equation* $x^2 + ky^2 = z^2$ *and the integral triangles with a cosine value of* $\frac{1}{n}$, Mathematics and Computer Education, Fall 2006